\documentclass[reqno,twoside, 10pt]{amsart}
\usepackage{amsmath}
\usepackage{color}
\usepackage{amsthm,amssymb,amsfonts}
\usepackage{varioref}
\usepackage[bookmarksnumbered,colorlinks,plainpages,backref]{hyperref}

\newtheorem{thm}{Theorem}

\newtheorem{prop}[thm]{Proposition}

\newcommand{\ue}{u_\varepsilon}

\newcommand{\lme}{{\lambda^\varepsilon}}

\renewcommand{\varrho}{\rho}
\newcommand{\e}{\varepsilon}

\newcommand{\lm}{\lambda}
\newcommand{\intl}{\int\limits}
\newcommand{\D}{\displaystyle}

\newcommand{\cL}{{\mathcal{L}}}

\newcommand{\cE}{{\mathcal{E}}}

\newcommand{\f}{\varphi}
\newcommand{\tg}{\mathop{\rm tg}\nolimits}
\def\D{\displaystyle}

\begin{document}
\author{Natalia Babych}
\title[High frequency asymptotics]{High frequency asymptotics of global vibrations in a problem with concentrated mass}
\date{}

\bigskip

\noindent
\small{%
%
N.Babych, \emph{High frequency asymptotics of global vibrations in a problem with con\-centra\-ted mass},
Math. Methods \& Phys - Mech. Fields {\bf 42}(3): 36 - 44, 1999.
\vglue1cm}

\begin{abstract}
We consider an elastic system containing a small region where the density is very much higher then elsewhere. Such system possesses two types of
eigenvibrations, which are local and global vibrations. Complete asymp\-totic
expansions of global eigenvibrations for ordinary differential opera\-tor of the fourth order are constructed using WKB~--~technique.

\bigskip


\medskip


\end{abstract}

\maketitle

\vspace{2mm}
{\bf Introduction.}
Heterogeneous systems give rise to new effects that do not reside in the separated system parts.
For instance, problems with local density perturbation are characterized by a presence of local and global vibrations.
For the first time the effects have been described by E.~Sanchez-Palencia~[6,7]. Vibrating systems with a local mass perturbation
are investigated starting from O.~A.~Oleinik~[3,4]. The effect of local vibrations has been studied for wide class
of the systems, for instance see~[1,2]. The global vibrations have been remained to be weakly analyzed so far.
A complete asymptotic expansions of the global eigenvibrations for one-dimensional system of the forth order with
locally perturbed density are constructed in this work.

\vspace{2mm}
{\bf 1. Problem statement.}\label{sec1}
Let a differential expression $L$ be given by 
$$
L v=(k_0(x)v'')''-(k_1(x)v')'+k_2(x)v,
$$
where the functions $k_0>0$, $k_1$, $k_2\ge 0$ are smooth at $[a,b]$. Denote by $[v]_{x=c}$
a jump of a function $v$ at point $c$. We investigate asymptotic behaviour as $\e\to0$
of eigenvalues $\lm_\e$ and eigenfunctions $u_\e$ of the problem
\begin{eqnarray}
 &Lu_\e-\lm_\e p(x)u_\e=0,\quad x\in (a,-\e)\cup(\e,b),  & \label{eq1}\\
 &Lu_\e-\lm_\e \e^{-m}q(x/\e)u_\e=0,\quad x\in (-\e,\e), & \label{eq2}\\
 &u_\e(a)=u'_\e(a)=0,\quad u_\e(b)=u'_\e(b)=0,           & \label{bc}\\
 &\left[u_\e\right]_{x=\pm\e}=[u'_\e]_{x=\pm\e}=[u''_\e]_{x=\pm\e}=
                   [u'''_\e]_{x=\pm\e}=0.  & \label{sp}
\end{eqnarray}
For each fixed $\e>0$ the problem possesses a countable set of eigenvalues.
Bihaviour as $\e\to0$ for each eigenvalue $\lme$ and corresponding eigenfunctions $\ue$
depending on a value of real parameter $m$ is investigated in~[1].
In the case $m>4$ problem~(\ref{eq1})~--~(\ref{sp}) possesses local eigenvibrations\footnote{We use a term ''eigenvibration'' 
to denote a pare of eigenvalue and corresponding eigenfunction} with
corresponding eigenvalues $\lm_n^\e=O(\e^{m-4})$, $\e\to 0$, and with eigenfunctions $u_n^\e$, which
localize itself in the interval of density perturbation $(-\e,\e)$, rapidly vanishing outside the interval.
Nevertheless, qualitative behaviour of the vibrating system is not yet described completely by local eigenvibrations.
Models with concentrated masses possess also global vibrations~[7]. As indicated below, the global
vibrations are supported by eigenfunction sequences with nontrivial limits $u_{n(\e)}^{\e}\to v_0$
for eigenvalues $\lm_{n(\e)}^\e\to \lm_0>0$ with $n(\e)\to \infty$.
Dependence $n(\e)$ is a discrete one that causes construction of asymptotics along certain sequences of the small parameter only.
A family of the sequences is bound by a deformation parameter, which is present in the lower terms of the constructed expansions.
Asymptotics depend on the value of $m$. We choose $m=8$ as a pattern case.

\vspace{2mm}
{\bf 2. Asymptotics of global vibrations: the leading terms.}\label{sec2}
We seek the asymptotic expansions of the eigenvalues $\lm_\e$ and the eigenfunctions $u_\e$
of problem~(\ref{eq1})~--~(\ref{sp}) in the form
\begin{eqnarray}
&\lm_\e\sim \D\sum\limits_{i=0}^{\infty}\e^i\lm_i,&\label{om}\\
&u_\e(x)\sim \D\sum\limits_{i=0}^{\infty}\e^iv_i(x),\quad x\in (a,-\e)\cup(\e,b),
\quad v_0\not\equiv 0.&\label{v}
\end{eqnarray}
In order to explore the function $\ue$ in the region $(-\e,\e)$, we consider
problem~(\ref{eq1})~--~(\ref{sp}) in variables $\xi=\e^{-1}x$.
Taking into account
$$
k_i(\e\xi)= \sum\limits_{j=0}^{\infty}(\e\xi)^j k_i^{(j)}(0)(j!)^{-1},
$$
we transform the differential expression $L$ to
$L_\e=\e^{-4} \D\sum\limits_{j=0}^{\infty}\e^{j}{\mathcal{ L}}(j)$,
where
$$
{ \mathcal{L}}(j)=\frac{k_0^{(j)}(0)}{j!}\frac{d^2}{d\xi^2}\xi^j\frac{d^2}{d\xi^2}-
\frac{k_1^{(j-2)}(0)}{(j-2)!}\frac{d}{d\xi}\xi^{j-2}\frac{d}{d\xi}+
\frac{k_2^{(j-4)}(0)}{(j-4)!}\xi^{j-4}.
$$
We use notation $k_n^{(i)}(x)$ for the $i$-th derivative of a function $k_n(x)$ in case if $i\in\mathbb{N}$,
$k_n^{(0)}(x)=k_n(x)$ and $k_n^{(j)}(x)\equiv 0$ if $j<0$.

Hence, the eigenfunction $\ue$, which corresponds to $\lme$, is a solution of the equation with
small parameter nearby the highest derivative
\begin{equation}
 \e^8L_\e U_\e-\lm_\e q(\xi)U_\e=0,\quad U_\e=U_\e(\xi),\quad\xi\in (-1,1).
\label{WKB}
\end{equation}
According to the method of WKB-approximations~[4], we seek the solution of \eqref{WKB}
as linear combination of the series
$e^{{\D\e^{-1}}S(\xi)} \D\sum\limits_{i=0}^{\infty}\e^ia_i(\xi)$.
Equality \eqref{WKB} is guarantied by the choice of functions $S$ and $a_i$.
In particular, the phase function $S$ is a solution of the eikonal equation
\begin{eqnarray}
k_0(0)S'^4-\lm_0 q(\xi)=0,\quad\xi\in(-1,1).\label{S}
\end{eqnarray}
Therefore we fix
$$
S(\xi)=(\lm_0k_0^{-1}(0))^{1/4} \intl_{-1}^\xi q^{1/4}(t)\,dt
$$
and, without loss of generality, we derive the eigenfunction $\ue$ in the form
\begin{equation}
u_\e(\e\xi)\stackrel{def}=U_\e(\xi)=\e^4 \D\sum\limits_{i=0}^{\infty}\e^i
\langle f_i(\xi),N(\xi,\e^{-1}S)\rangle,\quad\xi\in(-1,1),
\label{U}
\end{equation}
with $f_i$ being a vector-function with values in $\mathbb{R}^4$,
$\langle\cdot,\cdot\rangle$ being the standard scalar product in $\mathbb{R}^4$,
and the operator $N$ mapping $[-1,1]\times C^\infty [-1,1]$ in $\mathbb{R}^4$ according to
$$
 N(\xi,\tau)=(\cos\tau(\xi), \sin\tau(\xi),
  \exp(-\tau(\xi)+\tau(-1)), \exp(\tau(\xi)-\tau(1))\,).
$$
The choice of the multiplier $\e^4$ and the argument shifts in exponents of series~(\ref{U})
is related to a certain normalization of the eigenfunction $\ue$, which is described in Section 4 (all other eigenfunctions differ
by a constant multiplier).

We construct formal asymptotic expansions~(\ref{om}),~(\ref{v}),~(\ref{U}) satisfying all conditions of problem~(\ref{eq1})~--~(\ref{sp}).
In particular, equality~(\ref{eq1}) holds if
\begin{eqnarray}
&Lv_0-\lm_0 p(x)v_0 =0,\quad x\in(a,0)\cup(0,b),&\label{eqv0}\\
&Lv_i-\lm_0 p(x)v_i =p(x) \D\sum\limits_{j=1}^i\lm_jv_{i-j},\quad x\in(a,0)\cup(0,b).&
\label{eqv}
\end{eqnarray}
Let consider the action of the differential expression $L_\e$ on the function $U_\e$.
We note that
${\D\frac{d}{d\xi}}N(\xi,\tau)=\tau'(\xi)TN(\xi,\tau)$,
where the matrix
$$
T=\left(\begin{array}{cc}
            {T_1}&{0}\\
        {0}&{T_2}
        \end{array}
  \right),\quad
T_1=\left(\begin{array}{cr}
            {0}&{-1}\\
            {1}&{0}
        \end{array}
  \right),\quad
T_2=\left(\begin{array}{rc}
            {-1}&{0}\\
            {0}&{1}
        \end{array}
  \right)
$$
has the property $T^*=T^3$, and $T^4$ is a unit matrix.
Counting
$$
 {\frac{d}{d\xi}}\langle f_k(\xi),N(\xi,\e^{-1}S)\rangle=
 \Big\langle \Big(\e^{-1}S'(\xi)T^3+{\frac{d}{d\xi}}\Big)f_k(\xi),
N(\xi,\e^{-1}S)\Big\rangle,
$$
we get
$$
L_\e U_\e=
 \D\sum\limits_{i=0}^{\infty}\e^{i-8} \D\sum\limits_{j=0}^i\Big\langle
 \D\sum\limits_{l=0}^4{ \mathcal{L}}_l(j-l)f_{i-j}(\xi),N(\xi,\e^{-1}S)\Big\rangle,
$$
where ${ \mathcal{L}}_n(k)$ are the differential expressions of the $n$-th order with coefficients
depending on $S'$ and $\xi^{k}$. In particular,
$
 { \mathcal{L}}_0(i)=k_0^{(i)}(0)(i!)^{-1}S'^4\xi^i
$,
and
$$
 { \mathcal{L}}_1(i)=k_0^{(i)}(0)(i!)^{-1}
\Big (S'^3\xi^i{\D\frac{d}{d\xi}}+S'^2\xi^i{\D\frac{d}{d\xi}}S'+
 S'{\D\frac{d}{d\xi}}S'^2\xi^i
 +{\D\frac{d}{d\xi}}S'^3\xi^i\Big)T.
$$
We recall that the differential expressions ${ \mathcal{L}}_n(t)$, $n=0,\dots,4$,
are equal to zero one for $t<0$. Hence, equation~(\ref{eq2}) yields
$$
({ \mathcal{L}}_0(0)-\lm_0 q)f_i+({ \mathcal{L}}_0(1)+{ \mathcal{L}}_1(0)-\lm_1q)f_{i-1}=\chi_i,
$$
with $ \chi_i=- \D\sum\limits_{j=2}^i \D\sum\limits_{l=0}^4({
\mathcal{L}}_l(j-l)-\lm_jq)f_{i-j} $ for $i\ge 2$, and $\chi_t\equiv
0$ for $t<2$. The expression ${ \mathcal{L}}_0(0)-\lm_0 q$ is equal
to zero because it is the left-hand side of the eikonal
equation given by~(\ref{S}). Therefore the coefficients $f_i$ of
expansion~(\ref{U}) solve the system of the
first-order differential equations
\begin{equation}
({ \mathcal{L}}_0(1)+{ \mathcal{L}}_1(0)-\lm_1q)f_i=\chi_{i+1},\quad i=0,1,2,\dots\, .
\label{f}
\end{equation}
Boundary conditions~(\ref{bc}) give
\begin{equation}
v_i(a)=v'_i(a)=0,\quad v_i(b)=v'_i(b)=0,\quad i=0,1,2,\dots\,.
\label{bcv}
\end{equation}
Interfacial conditions~(\ref{sp}) applied to~(\ref{v}) and~(\ref{U}) yield
\begin{eqnarray}
\begin{array}{c}
  \D\sum\limits_{j=0}^i(\pm 1)^j(j!)^{-1}v_{i-j}^{(j)}(\pm 0)=
   \langle f_{i-4}(\pm1),N(\pm1,\e^{-1}S)\rangle, \\
  \D\sum\limits_{j=0}^i(\pm 1)^j(j!)^{-1}v_{i-j}^{(j+1)}(\pm 0)=
   \langle (S'T^3f_{i-2}+f'_{i-3})(\pm1),N(\pm1,\e^{-1}S)\rangle ,
\end{array} \label{sp1}\\[0.1cm]
\begin{array}{c}
  \D\sum\limits_{j=0}^i(\pm 1)^j(j!)^{-1}v_{i-j}^{(j+2)}(\pm 0)=
   \langle (S'^2T^2f_i+D_i)(\pm1),N(\pm1,\e^{-1}S)\rangle, \\
 F_i(\pm0)=\langle (S'^3Tf_i+E_i)(\pm1),N(\pm1,\e^{-1}S)\rangle,
\end{array} \label{sp2}
\end{eqnarray}
where we used notation 
$$
f_j\equiv 0\quad \text{for}\quad j<0;
$$
$$
D_i=2S'T^3f'_{i-1}+S''T^3f_{i-1}+f''_{i-2}
\quad \text{for}\quad i=0,1,\dots;
$$
$F_l$ take zero values for $l=0,1$, and
$$
F_i(\pm0)= \sum\limits_{j=0}^{i-2}(\pm1)^j(j!)^{-1}
v_{i-j-2}^{(j+3)}(\pm0) \ \
\mbox{for} \ \ i=2,3,\dots;
$$
and
$
E_i=\left(S'^2 {\D\frac{d}{d\xi}} + S' {\D\frac{d}{d\xi}} S'+
{\D\frac{d}{d\xi}} S'^2 \right)
T^2f_{j-1}+\left(S' {\D\frac{d^2}{d\xi^2}}+
{\D\frac{d}{d\xi}}S'{\D\frac{d}{d\xi}}+
{\D\frac{d^2}{d\xi^2}}S'
    \right)\times$\\
$T^3f_{j-2}+f'''_{j-3}
$
for $i=0,1,\dots$ .

According to conditions~(\ref{eqv0}),~(\ref{bcv}),~(\ref{sp1}) for $i=0$,
the leading terms $\lm_0$ and $v_0$ of expansions~(\ref{om}),~(\ref{v})
are correspondingly an eigenvalue and eigenfunction of the problem
\begin{equation}\label{v0}
\begin{array}{c}
  Lv_0-\lm_0 p(x)v_0=0,\quad x\in(a,0)\cup(0,b),\\
  v_0(a)=v'_0(a)=0,\quad v_0(0)=v'_0(0)=0,\quad v_0(b)=v'_0(b)=0.
\end{array}
\end{equation}
We restrict ourself to considering only simple eigenvalues of three-point Dirichlet problem~(\ref{v0}).
Nevertheless, there are situations when all the eigenvalues has multiplicity more then $1$ (for instance,
$a=-b$ and $p$ is an even function). We fix a simple eigenvalue $\lm_0$ of~(\ref{v0})
and the corresponding eigenfunction $v_0$ such that $v_0(x)\equiv 0$
for $x\in(0,b)$ and $\D\intl_a^0 pv_0^2\,dx=1$.

Equalities~(\ref{eqv}),~(\ref{bcv}),~(\ref{sp1}) for $i=1$ give
\begin{equation}\label{v1}
\begin{array}{c}
  Lv_1-\lm_0 p(x)v_1=\lm_1 p(x)v_0,\quad x\in(a,0)\cup(0,b),\\
  v_1(a)=v'_1(a)=0,\quad v_1(-0)=0,\quad v'_1(-0)=v''_0(-0),\\
  v_1(+0)=v'_0(+0)=0,\quad v_1(b)=v'_1(b)=0.
\end{array}
\end{equation}
Three-point problem~(\ref{v1}) has solution if and only if the parameter $\lm_1$
takes value
$$
 \lm_1=k_0(0)v''_0(-0)^2.
$$
We fix a solution $v_1$ of~(\ref{v1}) such that
$ \D\intl_a^0pv_1v_0\,dx=0$. One can note that $v_1(x)\equiv 0$ for $x\in(0,b)$.

Conditions~(\ref{f}),~(\ref{sp2}) for $i=0$ yield that the leading term $f_0$
of series~(\ref{U}) satisfies the problem
\begin{eqnarray}
  & f'_0 =A(\xi)f_0,\quad \xi\in(-1,1),& \label{A}\\
  & \langle T^2f_0(\pm1),N(\pm1,\e^{-1}S)\rangle=\sigma^{\pm},\quad
  \langle Tf_0(\pm1),N(\pm1,\e^{-1}S)\rangle=0,& \label{Abc}
\end{eqnarray}
with $\sigma^{+}=0$, $\sigma^{-}=S'(-1)^{-2}v''_0(-0)$.
The matrix of the linear homogeneous system of the first-order differential equations~(\ref{A}) is
$$
\begin{array}{c}
A=\left(\begin{array}{cc}
            {A_1}&{0}\\
        {0}&{A_2}
        \end{array}
  \right),\quad
A_1=\left(\begin{array}{rc}
            {\eta}&{\theta}\\
        {-\theta}&{\eta}
        \end{array}
  \right),\quad
A_2=\left(\begin{array}{cc}
            {\eta-\theta}&{0}\\
        {0}&{\eta+\theta}
        \end{array}
  \right),
\end{array}
$$
with
$$\eta=-3q'(\xi)(8q(\xi))^{-1},\
\theta= 1/4\,\,\left(\lm_0^{-3}k_0^{-5}(0)q(\xi)\right)^{1/4}
(\lm_1k_0(0)-\lm_0k'_0(0)\xi).$$

On the one hand, problem~(\ref{A}),~(\ref{Abc}) depends on the small parameter,
on another hand, it is a boundary value problem, and hence, it is an ill-posed one.
Both of the difficulties can be solved by exploring problem~(\ref{A}),~(\ref{Abc})
along with similar problems, which appear below, along discrete sequences of the small parameter
$\e_l\to 0$. Then "along" these sequences the problems have a unique solution, which
is independent of $\e_l$ up to the exponentially small terms. In particular, this property is consistant with
a discrete character of the high-frequency-vibrations effect, which is under consideration.

\begin{prop}\label{prop1}
Let $w:[-1,1]\to \mathbb{R}^4$ be a smooth vector-function, and $\sigma$ be a vector in $\mathbb{R}^4$.
There exists a small parameter sequence $\{\e_l\}_{l=1}^\infty$ such that the problem
\begin{eqnarray}
 &y'(\xi,\e)=A(\xi)\,y(\xi,\e)+w(\xi),\quad \xi\in(-1,1),& \label{y}\\[0.1cm]
&\begin{array}{cc}
 \langle T^2y(-1,\e),\,N(-1,\e^{-1}S)\rangle=\sigma_1,&
 \langle T y(-1,\e),\,N(-1,\e^{-1}S)\rangle=\sigma_2,\\
 \langle T^2y(1,\e),\,N(1,\e^{-1}S)\rangle=\sigma_3,&
 \langle T y(1,\e),\,N(1,\e^{-1}S)\rangle=\sigma_4,
\end{array}&\label{ybc}
\end{eqnarray}
has a unique solution $y(\cdot,\e_l)$ for each $l=1,2,\dots$.

The family of solutions $\{y(\cdot,\e_l)\}_{l=1}^\infty$ satisfies the estimate
$$
\| y(\cdot,\e_l)-y_*\|_{C^1}\le C e^{{\D -\e_p^{-1}}M}
$$
with $y_*$ being a smooth vector function on $[-1,1]$, and with positive constants $C$ and $M$ that
are independent of $\e$.
\end{prop}
\vspace{2mm}

P r o o f. \
System~(\ref{y}) has a fundamental matrix
$$
\begin{array}{c}
\Phi=q^{-3/8}
  \left(\begin{array}{cc}
            {\Phi_1}&{0}\\
        {0}&{\Phi_2}
        \end{array}
  \right),
\Phi_1=\left(\begin{array}{rc}
            {\cos\alpha}&{\sin\alpha}\\
        {-\sin\alpha}&{\cos\alpha}
        \end{array}
  \right),
\end{array}$$
$$ \begin{array}{c}
\Phi_2=\left(\begin{array}{cc}
            {e^{\alpha(-1)-\alpha}}&{0}\\
        {0}&{e^{\alpha-\alpha(1)}}
        \end{array}
  \right)
\end{array}
$$
with a function 
$$
 \alpha(\xi)={\frac{1}{4}}(\lm_0^{-3}k_0^{-5}(0))^{1/4} \intl_{-1}^{\xi}q^{1/4}(t)
 (\lm_1k_0(0)-\lm_0k'_0(0)t)dt
.$$
Then the general solution of~(\ref{y})  can be written as
$$
y(\xi)=\Phi(\xi)(\beta+h(\xi)),
$$
where $\beta$ is a constant vector,
and
$$
h(\xi)={\D\int}_{-1}^\xi\Phi^{-1}(t)w(t)\,dt.
$$
Let further 
$$
y(\xi,\e)=\Phi(\xi)(\beta_\e+h(\xi))
$$
be a solution of boundary value problem~(\ref{y}),~(\ref{ybc}).
It is easy to check that
$$
 \Phi^t(\xi)N(\xi,\e^{-1}S)=q^{-3/8}(\xi)N(\xi,\gamma_\e),\quad
 \mbox{ where }\, \gamma_\e(\xi)=\e^{-1}S(\xi)+\alpha(\xi),
$$
and $\Phi^t$ is transposed to $\Phi$. Therefore we get
$$
 \langle y(\xi,\e),N(\xi,\e^{-1}S)\rangle
  =q^{-3/8}(\xi)\langle \beta_\e+h(\xi),N(\xi,\gamma_\e)\rangle,
$$
that allows to rewrite boundary conditions~(\ref{ybc}) to the form
\begin{equation}\label{b}
\begin{array}{cl}
 \langle \beta_\e,\,T^2N(-1,\gamma_\e)\rangle &=m^{-}\sigma_1,\\
 \langle \beta_\e,\,T^3N(-1,\gamma_\e)\rangle &=m^{-}\sigma_2,\\
 \langle \beta_\e,\,T^2N(1,\gamma_\e)\rangle &
         =m^{+}\sigma_3-\langle h(1),\,T^2N(1,\gamma_\e))\rangle,\\
 \langle \beta_\e,\,T^3N(1,\gamma_\e)\rangle &
         =m^{+}\sigma_4-\langle h(1),\,T^3N(1,\gamma_\e))\rangle,
\end{array}
\end{equation}
where
$m^{\pm}=q^{3/8}(\pm1)$.
Here we have used the equalities
$T\Phi^t=\Phi^tT$ and $h(-1)=0$.
Therefore, the vector $\beta_\e$ has to be a solution to a linear algebraic system with matrix
$$
 G(\gamma_\e(1))=\left(\begin{array}{cccc}
           {-1}&{0}&{1}&{e^{-{\D\gamma_\e}(1)}}\\
               {0}&{-1}&{-1}&{e^{-{\D\gamma_\e}(1)}}\\
               {-\cos\gamma_\e(1)}&{-\sin\gamma_\e(1)}&{e^{-{\D\gamma_\e}(1)}}&{1}\\
               {\sin\gamma_\e(1)}&{-\cos\gamma_\e(1)}&{-e^{-{\D\gamma_\e}(1)}}&{1}
               \end{array}\right).
$$
Note that the determinant
$$
\det G(\gamma_\e(1))=-2\cos\gamma_\e(1)+2e^{-{\D\gamma_\e}(1)}
\bigl(2-e^{-{\D\gamma_\e}(1)}\cos\gamma_\e(1)\bigr)
$$
becomes zero on an infinitely small sequence $\e$ (but not for any $\e>0$
since $\gamma_\e(1)\to\infty$ as $\e\to0$).

We fix a number $\delta$ from the interval $[0,2\pi)$ such that $\delta$
is different from $\pi/2$ and $3\pi/2$.
We construct the sequence $\e_l$ using the set of conditions
$\gamma_{{\D\e_l}}(1)=\delta+2\pi l$ for  $l=1,2,\dots$~. In other words,
\begin{equation}\label{e}
 \e_l(\delta)=\frac{S(1)}{\delta+2\pi l-\alpha(1)}
\end{equation}
for all $l\ge l_0$, where $l_0$ is the smallest natural number such that the denominator
of fraction~(\ref{e}) is positive. We use notation $\gamma_l=\gamma_{{\D\e_l}}$.
Using $\gamma_l\to\infty$ as $\e_l\to0$, we can treat the matrix
$G(\gamma_l)$ as an exponentially small perturbation of matrix
$$
 G_\delta=\left(\begin{array}{cccc}
           {-1}&{0}&{1}&{0}\\
               {0}&{-1}&{-1}&{0}\\
               {-\cos\delta}&{-\sin\delta}&{0}&{1}\\
               {\sin\delta}&{-\cos\delta}&{0}&{1}
               \end{array}\right).
$$
In the same manner
$N(1,\gamma_l)=\left(\cos\delta, \sin\delta, e^{-{\D\gamma_l}(1)}, 1 \right)$,
then the right hand side of system~(\ref{b}) also is an exponentially small perturbation
of the vector
$$
g=(m^-\sigma_1,\, m^-\sigma_2,\, m^+\sigma_3-\langle h(1),\,T^2N_1\rangle ,
\, m^+\sigma_4-\langle h(1),\, T^3N_1\rangle ),
$$
where $N_1$ differs from $N(1,\gamma_l)$ only by the third coordinate being zero.

For the chosen value $\delta$ the matrix $G_\delta$ is a non-degenerate one. Then
according to the results of perturbation theory in a finite dimensional space, we have
$$
 \|\beta_{{\D\e_l}}-\beta_*\|_{\mathbb R^4}\le C e^{-{\D\gamma_l}(1)},
$$
where
$\beta_*$ is a solution to the system $G_\delta\beta=g$.
Let $y_*(\xi)=\Phi(\xi)(\beta_*+h(\xi))$, then
$$
\| y(\cdot,\e_l)-y_*\|_{C^1}\le \|\Phi\|_{C^1}
\|\beta_{\e_l}-\beta_*\|_{\mathbb R^4},
$$
with $\|\Phi\|_{C^1}=\max\limits_{\xi\in[-1,1]}(\|\Phi(\xi)\|+\|\Phi'(\xi)\|)$.
To finish the proof, it remains to note that $\gamma_\e(1)\ge M\e^{-1}$ 
with a certain positive constant $M$.
\hfill $\Box$

\vskip.3cm

\textbf{Remark.}
The function $y_*$ is referred to as the principal solution of problem~(\ref{y}),~(\ref{ybc})
since constructing asymptotic expansions in power scale of $\e$ we can neglect exponentially small terms.
Nevertheless, the choice of the sequence $\e_l$, and then the choice of the solution $y_*$,
is not unique since it depends on $\delta$.
As we mentioned at the beginning, the latter is connected to the fact that
the discrete approximation of global vibrations allows deformation, hence the approximation is not unique.
A presence of the deformation parameter $\delta$ in asymptotics corresponds to the problem content.

\vskip.3cm

We come back to exploring problem~(\ref{A}),~(\ref{Abc}), which is problem~(\ref{y}),~(\ref{ybc})
with the right-hand side $w=0$ and
$\sigma=(S'(-1)^{-2}v''_0(-0),0,0,0)$.
The system is homogeneous, and then $f_0(\xi)=\Phi(\xi)\beta_0$, where
$$
 \beta_0=1/2\,\,q^{3/8}(-1)S'(-1)^{-2}v''_0(-0)
 \left(\tg\delta-1,-\tg\delta-1,\tg\delta+1,-1/\cos\delta\right)
$$
is a solution of the corresponding linear system with matrix $G_\delta$.

Hence, we have found the terms
 $\lm_0$, $v_0$, $\lm_1$, $v_1$, $f_0$ of expansions~(\ref{om}),~(\ref{v}),~(\ref{U}).
We recall that $f_0$ depends on value $\delta$. The construction of complete asymptotic expansions
for solutions to problem~(\ref{eq1}) --~(\ref{sp})
is conducted on the sequence
${\cE}_{\delta}=\{\e_l(\delta)\}_{l_0}^{\infty}$
that is chosen accordingly to~(\ref{e}).

\vskip.3cm

\vspace{2mm}
{\bf 3. Complete asymptotics of global vibrations.}
\label{sec3}
Let find coefficients $\lm_i$, $v_i$, $f_{i-1}$
of asymptotic expansions~(\ref{om}),~(\ref{v}),~(\ref{U}) for $i\ge 2$.
Using conditions (\ref{eqv}),~(\ref{bcv}),~(\ref{sp1})
we construct a boundary value problem for the general term $v_i$ in the form
\begin{equation}\label{vi}
\begin{array}{c}
  Lv_i-\lm_0p(x)v_i=p(x) \D\sum\limits_{j=1}^i\lm_jv_{i-j},\quad x\in(a,0)\cup(0,b),\\
  v_i(a)=v'_i(a)=0,\quad  v_1(b)=v'_1(b)=0, \\
  v_i(\pm0)=V_i(\pm0),\quad v'_i(\pm0)=W_i(\pm0).
\end{array}
\end{equation}
The right-hand sides of the boundary conditions in problem~(\ref{vi})
with the precision up to exponentially small terms are

$$
 V_i(\pm0)=\langle q^{-3/8}\Phi^{-1}f_{i-4}(\pm1),N_{\pm1}\rangle-
  \sum\limits_{j=0}^{i-1}(\pm1)^j(j!)^{-1}v_{i-j}^{(j)}(\pm0)
,$$
and
$$
 W_i(\pm0)=\langle q^{-3/8}\Phi^{-1}
    (S'T^3f_{i-2}+f'_{i-3})(\pm1),N_{\pm1}\rangle-
     \sum\limits_{j=0}^{i-1}(\pm1)^j(j!)^{-1}v_{i-j}^{(j+1)}(\pm0)
,$$ where $N_{-1}=(1,0,1,0)$, and the vector $N_1$ was introduces
during the proof of Proposition \ref{prop1}. Since $\lm_0$ is an
eigenvalue of problem~(\ref{v0}), boundary value problem~(\ref{vi})
could have no solution. The necessary and sufficient condition for
existence of the solution is equality $
 \lm_i=(k_0v''_0W_i-(k_0v''_0)'V_i)(-0)
$.
We fix the solution $v_i$ to problem~(\ref{vi}) such that
$$
\D\intl_a^0v_iv_0=0.
$$

Then we find $f_{i-1}$ from the problem
\begin{equation}\label{Ai}
\begin{array}{l}
  f'_i =A(\xi)f_i+1/4\,\,k_0^{-1}(0)S'^{-3}T^3\chi_{i+1},
            \quad \xi\in(-1,1),\\
  \langle T^2f_i(\pm1),\,N(\pm1,\e^{-1}S)\rangle=
          -\langle(q^{-3/8}S'^{-2}\Phi^{-1}D_i)(\pm1),N_{\pm1}\rangle+ \\
  \qquad\qquad
  S'(\pm1)^{-2} \D\sum\limits_{j=0}^i(\pm1)^j(j!)^{-1}v_{i-j}^{j+2}(\pm0),\\
  \langle Tf_i(\pm1),N(\pm1,\e^{-1}S)\rangle=\\
  \qquad\qquad S'(\pm1)^{-3}
  \left(F_i(\pm0)-\langle(q^{-3/8}\Phi^{-1}E_i)(\pm1),N_{\pm1}\rangle\right).
\end{array}
\end{equation}
According to Proposition \ref{prop1}, the right-hand side of the
system is a smooth function. Then for the sequence of small
parameter ${\cE}_\delta$ there is a solution to problem~(\ref{Ai}).

Hence the algorithm of constructing coefficients of expansions~(\ref{om}),~(\ref{v}),~(\ref{U})
and a sequence of small parameter is
$$
 \lm_0\to v_0\to\lm_1\to v_1\to{\cE}_{\delta}\to f_0\to
 \dots\to \lm_i\to v_i\to f_{i-1}\dots\, .
$$
We recall that all the coefficients, except $\lm_0$, $v_0$ and $\lm_1$, $v_1$,
depend on the value of parameter $\delta$. In other words,
to each $\delta\in[0,2\pi)$, $\delta\not=\pi/2$, $3\pi/2$
we assign the corresponding asymptotic series that approximate the same
form of global vibrations~$v_0$.

\vspace{2mm}
{\bf 4. Justification of asymptotics.}\label{sec4}
The idea of justification is, using asymptotic series~(\ref{om}), to choose exactly the eigenfunction sequence of problem (\ref{eq1}) -- (\ref{sp})
such that the sequence models global vibrations and is approximated by series (\ref{v}),~(\ref{U}).

We fix numbers
$n\in\mathbb N$ and $\delta\in[0,2\pi)$ different from $\pi/2$ and $3\pi/2$.
We introduce a number sequence
$\{\lm_\e^{\{n\}}\}_{\D{\e\in{\cE}_{\delta}}}$
and a sequence of functions
$\{u_\e^{\{n\}}\}_{\D{\e\in{\cE}_{\delta}}}$.
Namely, for each element $\e$ of the discrete set $\cE_\delta$ we assume
$$
 \lm^{\{n\}}_\e=\lm_0+\e\lm_1+\dots\e^n\lm_n , \nonumber\\
$$
\begin{eqnarray*}
 &u^{\{n\}}_\e(x)=\left\{
       {\begin{array}{ll}
           v_0(x)+\e v_1(x)+\dots+\e^nv_n(x),& x\in (a,-\e)\cup(\e,b),\\
           \e^4\, \D\sum\limits_{i=0}^{n+2}\e^i\langle f_i(\e^{-1}x),
               N(\e^{-1}x,\e^{-1}S)\rangle ,& x\in (-\e,\e),
       \end{array}}
           \right. &\label{un}
\end{eqnarray*}
where the values $\lm_i$, functions $v_i$, vectors $f_i$ and
the set ${\cE}_{\delta}$ are constructed in Sections 2 and 3. 
\vspace{2mm}

\begin{prop}\label{prop2}
There exists an eigenvalue sequence $\{\lm^\e\}_{\D{\e\in{\cE}_{\delta}}}$ 
such that
$$
|\lm^{\{n\}}_\e-\lm^{\e}|\leq C_n\e^{n+1},\quad \e\in{\cE}_{\delta}.
$$
\end{prop}
\vspace{2mm}

P r o o f. \ Denote by $\rho_\e$ the density of the original
problem, which equals $\e^{-8}q(x/\e)$ on the interval $(-\e,\e)$,
and is $p(x)$ outside the interval. Problem~(\ref{eq1})
--~(\ref{sp}) considered in the weight space $L_2(\rho_\e,(a,b))$ is
equivalent to the problem for the self-adjoint operator
${\cL}_\e=\rho_\e^{-1}L$:
$$
  {\cL}_\e u_\e-\lm_\e u_\e=0.
$$
The domain of the operator ${\cL}_\e$ is
$$
 D({\cL}_\e)=\{\f\in H^4(a,b): \f(a)=\f'(a)=0,
\ \
 \f(b)=\f'(b)=0\}.
$$
The spectrum of the operator ${\cL}_\e$ is discrete [1].

The function $u^{\{n\}}_\e$ is not an element of the space $D({\cL}_\e)$
since it has discontinuities at the points $x=\pm\e$.
Nevertheless, we state existence of the function $\f^{\{n\}}_\e$
such that
$u^{\{n\}}_\e+\f^{\{n\}}_\e \in D({\cL}_\e)$,
moreover $\f^{\{n\}}_\e$ is equal to zero on the interval $(-\e,\e)$,
and
\begin{equation}\label{fi}
  {\D\max_{x\in (a,b)}}
   \D\sum\limits_{i=0}^4\left|{\D\frac{d^i}{dx^i}}\f^{\{n\}}_\e\right|\le C_n\e^{n+1}.
\end{equation}
We introduce function
$$V^{\{n\}}_\e=\kappa_\e(u^{\{n\}}_\e+\f^{\{n\}}_\e),$$
where the factor $\kappa_\e$ is chosen accordingly to
$$\|V^{\{n\}}_\e\|_{{L_2(\rho_\e,(a,b))}}=1.$$
We remark that the value $\kappa_\e$ is separated from zero by
a positive constant that does not depend on $\e$.

Taking into account~(\ref{fi}), we can show that
$$
 \left\| {\cL}_\e V_\e^{\{n\}}-\lm_\e^{\{n\}}
  V_\e^{\{n\}}\right\|_{\D L_2(\rho_\e,(a,b))}
 \leq K_n\e^{n+1}, \quad \e\in{\cE}_{\delta},
$$
where the constant $K_n$ is independent of $\e$.
Therefore there exists an eigenvalue $\lm^\e$ of the operator ${\cL}_{\e}$ such that
$$
 \left|\lm^\e-\lm_\e^{\{n\}}\right|\leq K_n\e^{n+1}. 
 \qquad \Box 
$$

\vskip.3cm

Using asymptotics of each eigenvalue $\lm^\e_n$ we can prove that
in a certain $d\e^4$-vicinity of points $\lm_\e$, which were chosen accordingly to Proposal \ref{prop2},
we don't have other eigenvalues of the initial problem.
As consequence, for $n\ge 5$ we have only one eigenvalue $\lm_\e$
that satisfies Proposition~\ref{prop2}.
We choose the sequence of exactly these $\{\lm_\e\}_{\D{\e\in{\cE}_\delta}}$
and the corresponding normalized eigenfunctions $\{u_\e\}_{\D{\e\in{\cE}_\delta}}$.
Then
$$
 \|u_\e-\kappa_\e v_\e^{\{n\}}\|_{\D L_2(\rho_\e,(a,b))}\le K_nd^{-1}\e^{n-3},\quad
 \|u_\e\|_{\D L_2(\rho_\e,(a,b))}=1.
$$
This inequality complete the proof of

\vspace{5mm}

\textbf{Theorem.}
\emph{On the sequence ${\cE}_\delta$ the chosen functions $u_\e$ satisfy the estimates
$$
 \left\|u_\e(x)-\kappa_\e \D\sum\limits_{k=0}^n\e^kv_k(x)\right\|_{\D L_2(\Omega_\e)}\le
 C_n\e^{n+1},
$$
where region $\Omega_\e$ states for $(a,-\e)$ or $(\e,b)$, and the estimate
$$
 \left\|u_\e(\e\xi)-\e^4\kappa_\e \D\sum\limits_{k=0}^n\e^k
 \langle f_k(\xi),\,N(\xi,\e^{-1}S)\rangle\right\|_{\D L_2(-1,1)}\le  C_n\e^{n+1}
$$
for $n=0,1,\dots$ .
The coefficients $v_k$, $f_k$ are the same as at the beginning of Section 4. 
The normalizing multiplier $\kappa_\e$ tends to $1$ as $\e\to0$.
}


\vspace{20mm}

\vbox{\noindent\small
Lviv National University 09.12.98 }

\end{document}